\documentclass{amsart}

\usepackage{amssymb,amsmath,epsfig}

\newtheorem{thm}{Theorem}
\newtheorem{lem}[thm]{Lemma}
\newtheorem{prop}[thm]{Proposition}
\theoremstyle{definition}

\theoremstyle{remark}
\newtheorem*{rem}{Remark}

\newcommand{\Z}{{\mathbf{Z}}}
\newcommand{\R}{{\mathbf{R}}}
\newcommand{\C}{{\mathbf{C}}}

\newcommand{\Rat}{{\bf Rat}}
\newcommand{\Rt}{{\bf \overline{Rat}}}

\newcommand{\CPN}[1]{{\bf CP}^{#1}}

\font\cy=wncyr10

\title{Spaces of Rational Maps and 
the Stone-Weierstrass Theorem}

\author{Jacob Mostovoy}

\address{Instituto de Matem\' aticas (Unidad Cuernavaca)
Universidad Nacional Aut\' onoma de M\' exico, Av.\ Universidad s/n, 
col. Lomas de Chamilpa, C.P. 62210, Cuernavaca, Morelos, M\'exico.}

\begin{document}

\begin{abstract} It is shown that Segal's theorem on the spaces of rational
maps from $\CPN{1}$ to $\CPN{n}$ can be extended to the spaces of continuous
rational maps from $\CPN{m}$ to $\CPN{n}$ for any $m\leqslant n$. The tools are
the Stone-Weierstra\ss\ theorem and Vassiliev's machinery of
simplicial resolutions.
\end{abstract}



\maketitle

\section{Introduction.}

Let $X$ and $Y$ be two topological spaces with some additional structure.
Maps from $X$ to $Y$ that preserve the structure form a subspace in the
space of all continuous maps from $X$ to $Y$. One may ask whether
this subspace gives some kind of a ``topological approximation''
to the space of all continuous maps.

For example,  if $X$ is a Stein manifold and $Y$ is a complex manifold with
a dominating spray (homogeneous spaces and complements in $\C^N$ of algebraic
subvarieties of codimension at least two satisfy this condition),
the inclusion of the space of holomorphic maps from $X$ to $Y$ into
the space of all continuous maps from $X$ to $Y$ is a weak homotopy
equivalence. This statement is a corollary of the Oka-Grauert principle
proved by Gromov \cite{Gromov}.

Another situation which has been extensively studied is the case when
$X$ and $Y$ are compact complex manifolds and $X$ has dimension one.
Segal in \cite{Segal} proved, among other things, that the inclusion of the
space of (based) rational maps of degree $d$ from $\CPN{1}$ to $\CPN{n}$
into the two-fold loop space of $\CPN{n}$ is a homotopy equivalence up to
dimension $(2n-1)d$. The results of \cite{Segal} were later extended
by various authors: see, for example
\cite{Boyer+},\cite{Guest-A},\cite{Guest-B},\cite{JM},\cite{Kallel+}. 
In all these generalizations
$X$ is a curve. Apparently, the only published works that address the case
$\dim{X}>1$ are the paper \cite{Havlicek} which deals with the space of
holomorphic maps from $\CPN{1}\times\CPN{1}$ to the Grassmannian of
two-dimensional planes in $\C^{N}$, and the work \cite{KYlin}
describing spaces of linear maps between complex projective spaces.

Segal conjectured in \cite{Segal} that the stability theorem
he proved for spaces of rational curves can be extended
to the spaces of maps from $\CPN{m}$ to $\CPN{n}$ for $m>1$.
Such a generalization will be discussed in this paper.

\medskip

A  continuous rational map from $\CPN{m}$ to $\CPN{n}$
can be given by a collection of $n+1$ complex homogeneous
polynomials of the same degree and with no common zero. Up to a constant,
such a representation  of a map by polynomials is unique. If $m>n$, all
continuous rational maps are constant, so we shall always assume
that $m\leqslant n$. By the degree of a rational map we shall mean
the degree of the polynomials that define it. More generally,
a continuous map between complex projective spaces
induces multiplication by an integer in the second cohomology.
We shall call this integer the  degree of the corresponding
map\footnote{this corresponds to the algebraic rather than
topological degree.}.
In what follows, all rational maps will be continuous.

Let $\Rat_{f}^{m,n}$ be the set of all rational maps of degree $d$
from $\CPN{m}$ to $\CPN{n}$ that restrict to a given map
\[ f:\CPN{m-1}\to\CPN{n}\]
on a fixed hyperplane $\CPN{m-1}\subset\CPN{m}$.
Fixing the coefficients
of the polynomials that define $f$ one obtains a bijection between
$\Rat_{f}^{m,n}$ and a subset of a complex affine space: a
rational map under this bijection is sent to the set of the coefficients of
the polynomials defining it. Thus $\Rat_{f}^{m,n}$ acquires a topology
from the complex affine space.

Consider the space of all continuous maps from $\CPN{m}$ to $\CPN{n}$ that
restrict to a given map $f$ on $\CPN{m-1}$. This space is
homotopy equivalent to the $2m$-fold loop space of $\CPN{n}$
and in what follows will be denoted by $\Omega^{2m}\CPN{n}$.

\begin{thm}\label{the}
The map
\[\Rat_{f}^{m,n}\to\Omega^{2m}\CPN{n}\]
given by the inclusion of rational maps into the space of all continuous maps
from $\CPN{m}$ to $\CPN{n}$ that restrict to a given map $f$ of degree $d$
on a fixed hyperplane, induces an isomorphism in homology 
in all dimensions smaller than 
\[ (2n-2m+1)\left(\left[\frac{d+1}{2}\right]+1\right).\]
If $m<n$ it is also induces isomorphisms of homotopy groups in these
dimensions. Here $[x]$ denotes the integer part of $x$.
\end{thm}

Theorem~\ref{the} implies a similar result for the spaces of ``free" maps.
Denote by $\text{\sl Map}_{d}(m,n)$ the space of all continuous maps 
from $\CPN{m}$ to $\CPN{n}$ 
of degree $d$, and let $\Rat_{(d)}^{m,n}$ be its subspace formed by all 
rational maps. 

\begin{thm} \label{free} In the statement of Theorem~\ref{the} the spaces 
$\Rat_{f}^{m,n}$ and $\Omega^{2m}\CPN{n}$ can be replaced by $\Rat_{(d)}^{m,n}$
and $\text{\sl Map}_{d}(m,n)$ respectively.
\end{thm}   

The proof of Theorem~\ref{the}
consists of constructing a sequence of topological spaces
and ``stabilization'' maps between them which starts with $\Rat_{f}^{m,n}$
and has the loop space $\Omega^{2m}\CPN{n}$ as its colimit. We shall then 
see that all
stabilization maps in this sequence induce isomorphisms in homology 
up to a certain dimension. This, together with the fact that for $m<n$ the
spaces  $\Rat_{f}^{m,n}$ are simply-connected, will imply Theorem~\ref{the}.

In the next section the stabilization maps are defined and it is shown
that the ``stabilized'' space of rational maps is homotopy equivalent to
$\Omega^{2m}\CPN{n}$. The trick behind our construction of the
stabilization maps is to place the problem in a real rather than
complex context and use the Stone-Weierstrass Theorem.
Section~\ref{SR} is a brief review of simplical resolutions.
These are then used in section~\ref{VSS} to construct
and describe the Vassiliev spectral sequences converging to the homology of
the spaces $\Rat_{f}^{m,n}$. We follow Vassiliev's methods as described in
\cite{V1,V2}, the only minor novelty being our use of degenerate
simplicial resolutions. The case of the free maps is treated in 
section~\ref{F}. Some comments are given in the last section.


\section{Stabilization maps and the Stone-Weierstrass Theorem.}

There are two ways of defining stabilization maps for the spaces of rational
maps. One of these constructions generalises Segal's stabilization maps;
we are not going to discuss it here. The stabilization we shall use can
be described as ``real stabilization''. In order to define it we need to
consider a wider class of maps between projective spaces.

Define a $(p,q)$-polynomial to be a homogeneous polynomial in
``holomorphic'' and ``anti-holomorphic'' variables which is
of degree $p$ in the holomorphic and of degree $q$ in the anti-holomorphic
variables. By a $(p,q)$-map we shall mean a map from $\CPN{m}$ to $\CPN{n}$
given by a collection of $n+1$ $(p,q)$-polynomials with no common zero. Here
the holomorphic variables are the homogeneous coordinates $z_i$ in $\CPN{m}$
and the anti-holomorphic variables are their conjugates.
In this terminology $(d,0)$-maps are precisely the rational maps of degree $d$.
In general, two representations of a $(p,q)$-map  by collections
of $(p,q)$-polynomials need not coincide up to a constant but
rather up to multiplication by a positive function. Thus a $(p,q)$-map
can be also thought of as a $(p+1,q+1)$-map:  just multiply each of the
polynomials defining the map by $z_0\bar{z}_0+\ldots +z_m\bar{z}_m$.
The degree of a $(p,q)$-map is readily seen to be equal to $p-q$.

Let $\Rat_{f}(p,q)$ be the space of $(p,q)$-maps from $\CPN{m}$ to
$\CPN{n}$ that restrict to a given map $f$ of degree $p-q$ on a hyperplane
$\CPN{m-1}\subset\CPN{m}$. In what follows we shall always assume that this
hyperplane is given by the equation $z_m=0$.
The topology on
$\Rat_{f}(p,q)$ is the topology of a subset of the
space of all continuous maps. Clearly, $\Rat_{f}(\deg{f},0)=\Rat_{f}^{m,n}$.
Considering $(p,q)$-maps as $(p+1,q+1)$-maps one obtains the inclusion
\[\Rat_f(p,q)\hookrightarrow\Rat_f(p+1,q+1).\]
We define $\Rat_f(d+\infty,\infty)$ to be the union of the spaces
$\Rat_f(d+k,k)$ for all $k\geqslant 0$.

The space of all continuous maps from $\CPN{m}$ to $\CPN{n}$
that restrict to $f$ on a fixed hyperplane is homotopy equivalent to the space
$\Omega^{2m}\CPN{n}$. Indeed, $\CPN{m}$ can be obtained from a $2m$-dimensional
disk $D^{2m}$ by making identifications on the boundary of $D^{2m}$; namely,
by collapsing the fibres of the Hopf map $\partial D^{2m}\to \CPN{m-1}$.
Hence, a map from $\CPN{m}$ to $\CPN{n}$ can be thought of as defined on $D^{2m}$;
the space of all maps from $D^{2m}\to \CPN{n}$ which restrict to the same map on 
$\partial D^{2m}$ is readily seen to be homotopy equivalent to $\Omega^{2m}\CPN{n}$.

\begin{prop}\label{th:real}
The natural inclusion of $\Rat_f(d+\infty,\infty)$ into the space \linebreak
 $\Omega^{2m}\CPN{m}$ of all continuous maps from $\CPN{m}$ to $\CPN{n}$
that restrict to $f$ on a fixed hyperplane, is a homotopy equivalence.
\end{prop}

Proposition~\ref{th:real} is a consequence of the following statement:
\begin{lem}\label{th:swt}
Let $X$ be a finite CW-complex. Any continuous map
\[ F:\CPN{m}\times X\to\CPN{n}\]
can be uniformly approximated with respect to the Fubini-Study metric on
$\CPN{n}$ by maps whose restriction to $\CPN{m}\times \{ x\}$ for any
$x \in X$, is a $(p,q)$-map for some $p,q$. Moreover, if the restriction
of $F$ to $\CPN{m-1}\times \{ x\}$ is a $(p,q)$-map for
all $x \in X$ , the approximating maps can be chosen so as
to coincide with $F$ on $\CPN{m-1}\times X$.
\end{lem}

%
\begin{proof}[Proof of Proposition~\ref{th:real}.]
For any compact riemannian manifold $M$ and any space $Y$ there
exists $\epsilon >0$ such that any two maps $Y\to M$ that are uniformly
$\epsilon$-close, are homotopic.
Thus, Lemma~\ref{th:swt} implies that the map of the homotopy groups
\[\pi_k\Rat_f(d+\infty,\infty)\to\pi_k\Omega^{2m}\CPN{n} \]
is surjective for all $k$. Indeed, setting $X=S^k$ we get that
any element of $\pi_k\Omega^{2m}\CPN{n}$ can be approximated by a homotopy
class in the space of $(p,q)$-maps. On the other hand, this homomorphism
is also injective, as any homotopy that goes through continuous maps can
be approximated by a homotopy through $(p,q)$-maps
(set $X=S^k\times [0,1]$). Both spaces of maps have homotopy types of
CW-complexes, so by the Whitehead Theorem they are homotopy equivalent.
%
\end{proof}

Lemma~\ref{th:swt} is a corollary of the Stone-Weierstra\ss\
Theorem for vector bundles:
\begin{thm}\label{th:SW}
Let $E$ be a locally trivial real vector bundle over a compact space $Y$,
$s_\alpha: Y\to E$ - a set of its sections, and let $A$ be a subalgebra
of the $\R$-algebra $C(Y)$ of continuous real-valued functions on $Y$.
Suppose that
\begin{itemize}
\item
the subalgebra $A$ separates points of $Y$, that is, for any
pair $x,y\in Y$ there exists $h\in A$ such that $h(x)\neq h(y)$;
\item
for any $y\in Y$ there exists $h\in A$ such that $h(y)\neq 0$;
\item for any $y\in Y$ the fibre of $E$ over $x$ is spanned by the $s_\alpha(y)$.
\end{itemize}
Then the $A$-module generated by the $s_\alpha$ is dense in the space of all
continuous sections of $E$.
\end{thm}

The above version of the Stone-Weierstra\ss\  Theorem follows from
the usual\footnote{see \cite{Rudin}.}  Stone-Weierstra\ss\  Theorem and the
existence of a partition of unity subordinate to the open cover of $Y$
trivialising $E$.

%
\begin{proof}[Proof of Lemma~\ref{th:swt}.]
Any continuous map $F:\CPN{m}\times X\to\CPN{n}$ can be given by
a collection of $n+1$ sections of some line bundle $E_F$ over 
$\CPN{m}\times X$. Namely, $E_F$ is the pullback of the tautological line 
bundle on $\CPN{n}$ with respect to $F$.

Choose an open cover of $X$ by contractible sets $U_{\beta}$ and 
a partition of unity $\rho_{\beta}$ subordinate to $U_{\beta}$, with $\beta$ 
belonging to some finite index set. Let $d$ be the degree of the 
restriction of $F$ to $\CPN{m}\times\{x\}$ for $x\in X$. The restriction
of $E_F$ to each subspace of the form $\CPN{m}\times U_{\beta}$ is 
isomorphic to the pullback of the $d$th power of the tautological line 
bundle on $\CPN{m}$; we shall assume that an explicit identification of 
these two line bundles is chosen for each $\beta$.

For each $\beta$ and each $i$ with $0\leqslant i \leqslant m$
let $s_{i,\beta}$ be the section of $E_F$ equal to $z_i^d\rho_{\beta}$ 
over $\CPN{m}\times U_{\beta}$ and trivial outside 
$\CPN{m}\times U_{\beta}$. Denote by $\sigma$ the set that consists of
the $s_{i,\beta}$ together with all the sections of the form 
$\sqrt{-1}s_{i,\beta}$ .

Let $A_0$ be the algebra of functions on $\CPN{m}$ 
generated by
\[\frac{z_i\bar{z}_i}{z_0\bar{z}_0+\ldots+z_m\bar{z}_m},\]
for $0\leqslant i\leqslant m$.
The first part of Lemma~\ref{th:swt} is then recovered
from Theorem~\ref{th:SW} applied to $E_F$:
the set of sections $s_\alpha$ is taken to coincide with $\sigma$
and $A$ is taken to be the algebra of all continuous functions on
$\CPN{m}\times X$ whose restriction to $\CPN{m}\times \{ x \}$ for any
$x\in X$ belongs to $A_0$.

In order to verify the second part of the lemma we shall show that
a sufficiently good approximation to $F$ can be modified so as to coincide
with $F$ on $\CPN{m-1}\times X$.

Let $S$ and $P$ be two sections of $E_F$ whose restrictions
to $\CPN{m-1}\times \{ x\}$ and $\CPN{m}\times \{ x\}$,
respectively, are $(p,q)$-polynomials for all $x\in X$.
Denote by $s$ and $p$ the restrictions of $S$ and $P$ to
$\CPN{m-1}\times X$. One can think of $s$ and $p$
as defined on all $\CPN{m}\times X$ and being independent of
$z_{m}$ and $\bar{z}_m$.

Consider a family $P_t$ of sections of $E_F$ with $t\in [0,1]$:
\[P_t=P+t(s-p).\]
For any $t$ the restriction of $P_t$ to any $\CPN{m}\times \{ x\}$
is a $(p,q)$-polynomial. Clearly, $P_0=P$ and the restriction of $P_1$
to any $\CPN{m-1}\times \{ x\}$ coincides with that of $S$.
Moreover, if $|P-S|<\epsilon$ then $|p-s|<\epsilon$ and, hence,
$|P_t-S|<2\epsilon$ in the standard metric that the fibres of $E_F$
inherit from the tautological line bundle on $\CPN{n}$.

Now, let $F_j$, $0\leqslant j\leqslant n$, be a set of  sections of $E_F$
that define the map $F$. Choose $\epsilon>0$ and find $\delta$ such that
any map $G$ from $\CPN{m}\times X$ to $\CPN{n}$ given by a set of
sections $G_j$ of $E_F$ with  $|F_j-G_j|<\delta$, is uniformly
$\epsilon$-close to $F$. Suppose that we have found $\delta/2$-approximations
$F'_j$ by families of $(p,q)$-polynomials for each of the $F_j$.
Applying the above construction with $P=F'_j$ and $S=F_j$ we
obtain an $\epsilon$-approximation to $F$ which coincides with $F$ on
$\CPN{m-1}\times X$ and whose restriction to any $\CPN{m}\times \{ x\}$
is a $(p,q)$-map.
%
\end{proof}


\section{Simplicial resolutions.}\label{SR}
Let $h:X\to Y$ be a finite-to-one surjective map of topological spaces
and let $\boldsymbol{i}$ be an embedding of $X$ into $\R^N$ for some $N$.
A {\em simplicial resolution} associated to the map $h$ with respect to
the embedding $\boldsymbol{i}$, is a subspace $X^{\Delta}$ of
$\R^N\times Y$ together with the projection map $h^{\Delta}:X^{\Delta} \to Y$.
The points of $X^{\Delta}$
are pairs $(t,y)$ with $y\in Y$ and $t$ belonging to the convex hull
of the set ${\boldsymbol{i}}\circ h^{-1}(y)$ in $\R^N$. The
space  $X$ is a subspace of $X^{\Delta}$; the restriction of $h^{\Delta}$
to $X$ coincides with the original map $f$. We say that a simplicial
resolution is {\em non-degenerate} if for each $y\in Y$ any $k$ points of the
set ${\boldsymbol{i}}\circ h^{-1}(y)$ span a $(k-1)$-dimensional affine
subspace of $\R^N$. Sometimes we shall use the term ``simplicial resolution''
for the space $X^{\Delta}$, this should not lead to confusion.

The fibres of the projection map $h^{\Delta}$ are contractible,
being convex polyhedra. We shall need simplicial resolutions in the situation
when $X$ and $Y$ are closed semialgebraic subsets of $\R^N$ and
$h$ and $\boldsymbol{i}$ are polynomial maps; from now on we shall assume
that this is the case. Under these circumstances $h^{\Delta}$
will always be a homotopy equivalence; hence, the problem of computing the
homotopy type of $Y$ is equivalent to the same problem for $X^{\Delta}$.

It is clear that any two non-degenerate simplicial resolutions associated to
the same map but with respect to different embeddings of $X$ into
$\R^N$ are homeomorphic over $Y$. This statement is generally false
without the non-degeneracy assumption.

There is an increasing filtration
\[ X_1\subset X_2\subset\ldots\subset X^{\Delta} \]
on any simplicial resolution associated to a map $h:X\to Y$.
Assume first that $h^{\Delta}:X^{\Delta} \to Y$ is non-degenerate. Then
for any $y\in Y$ its inverse image $(h^{\Delta})^{-1}(y)$
is a simplex; the subspace $X_k\subset X^{\Delta} $ is then defined
as the union of the $(k-1)$-skeleta of these simplices over all $y\in Y$.
In particular, $X_1=X$.

Now let $\widetilde{X}^{\Delta},X^{\Delta} \to Y$ be two simplicial
resolutions associated to the same map $h:X\to Y$
and suppose that $\widetilde{X}^{\Delta}\to Y$ is non-degenerate.
There exists the unique map
\[\pi:\widetilde{X}^{\Delta}\to{X}^{\Delta}\]
over $Y$ which commutes with the inclusions of $X$ into
$\widetilde{X}^{\Delta}$ and
${X}^{\Delta}$ and which is affine over all points of $Y$.
The increasing filtration on  ${X}^{\Delta}$ is defined
as the image of the increasing filtration on $\widetilde{X}^{\Delta}$
under $\pi$. This definition does not depend on the choice of a
non-degenerate resolution. Indeed, if both
$\widetilde{X}^{\Delta},{X}^{\Delta} \to Y$ are non-degenerate, the
canonical map from  $\widetilde{X}^{\Delta}$ to ${X}^{\Delta}$
not only is a homeomorphism, but also respects the increasing filtration.

The same construction can be carried out if the map $h:X\to Y$ is
not finite-to-one. However, in this situation every simplicial resolution,
as defined above, is necessarily degenerate. We define a non-degenerate
resolution as follows.

Let $\boldsymbol{i}_k:X\to\R^{N_k}$ be an embedding such that for any
$2k$ distinct points in $X$ their images in $\R^{N_k}$ span an affine
subspace of dimension $2k-1$. Then for all $j\leqslant k$,
for any $y\in Y$ and for any $j$ points in the set
$\boldsymbol{i}_k\circ h^{-1}(y)$
their convex hull is an $(j-1)$-simplex; moreover, the $(j-1)$-simplices
corresponding to disjoint sets of points in $\boldsymbol{i}_k\circ h^{-1}(y)$
are disjoint.
Now, define $X_k$ to be the union of all the convex hulls in $\R^{N_k}$
of the subsets of cardinality at most $k$ of
$\boldsymbol{i}_k\circ h^{-1}(y)$ for all $y\in Y$. There exists a natural
extension of $h$ to $X_k$.

For all $i<j$ the space $X_i$ can be considered as a subspace of $X_j$.
The non-degenerate simplicial resolution
$(\widetilde{X}^{\Delta},\widetilde{h}^{\Delta})$
associated to $h$ is then defined as the union of all the $X_i$ together
with the extension of $h$ to $\cup X_i$; the spaces $X_i$ form the
increasing filtration on $\widetilde{X}^{\Delta}$.

As before, the non-degenerate simplicial resolution does not depend on the
embeddings $\boldsymbol{i}_k$ and it is universal in the sense that for
any other simplicial resolution $X^{\Delta}$ there is a unique map
$\widetilde{X}^{\Delta}\to {X}^{\Delta}$ over $Y$ which is identity on
$X\in\widetilde{X}^{\Delta},X^{\Delta}$ and which is affine on the fibres.
Thus we have an increasing filtration on any simplicial resolution
of $h$ even in the case when $h$ is not finite-to-one.

Finally, let us state one obvious but important property of non-degenerate
simplicial resolutions. Assume there is a commutative square
\[ \begin{array}{ccc}
X            & \to & X'\\
\ \ \ \downarrow h &     &\ \ \ \downarrow h' \\
Y            & \to & Y'
\end{array} \]
with $h$ and $h'$ surjective. Then there is an induced filtration-preserving
map of non-degenerate simplicial resolutions associated to $h$ and $h'$.


\section{The Vassiliev spectral sequence.}\label{VSS}

The space  $\Rat_f(p,q)$ is defined as a subspace of $\Omega^{2m}\CPN{n}$.
In order to apply the machinery of simplicial resolutions we shall
replace it with a homotopy equivalent space $\Rt_f(p,q)$  defined as a
complement to an algebraic subvariety of an affine space.

Let us fix the coefficients of the $(p,q)$-polynomials $f_i$
defining the $(p,q)$-map $f:\CPN{m-1}\to\CPN{n}$.
Denote by
\begin{itemize}
\item[]$W^{i}_{p,q}$ - the complex affine space of all $(p,q)$-polynomials
that restrict to $f_i$ on the hyperplane $z_m=0$;
\item[]$W_{p,q}$ - the Cartesian product of the $W^{i}_{p,q}$ for
$0\leqslant i \leqslant n$;
\item[]$N_{p,q}$ - the complex dimension of $W_{p,q}$.
\end{itemize}

Define $\Rt_f(p,q)\subset W_{p,q}$ to
be the space of all $(n+1)$-tuples of $(p,q)$-polynomials which have no
common zero. The natural map
\[\Rt_f(p,q)\to\Rat_f(p,q)\]
is not a homeomorphism. However, it is obviously onto and its fibres are
convex and, hence, contractible. It is easy to see that the above map is,
in fact, a homotopy equivalence. The space $\Rt_f(p,q)$
is the complement of a discriminant $\Sigma$ in $W_{p,q}$, which consists
of $(n+1)$-tuples of $(p,q)$-polynomials that all have a common zero.
It has complex codimension $n-m+1$ in $W_{p,q}$ so the spaces
$\Rat_{f}^{m,n}$ are simply-connected if $m<n$.

In what follows the notation $\C^m$ will be used for
the affine chart $z_m=1$ in $\CPN{m}$. Let $V$ be the complex
vector space of dimension $\binom{p+m}{m}\binom{q+m}{m}$
spanned by all monomials in $z_i$ and $\bar{z}_i$
(with $0\leqslant i<m$) of degree at most $p$ in $z_i$ and
at most $q$ in $\bar{z}_i$. We shall denote by $v_{p,q}$
the Veronese-type embedding of $\C^m$ into $V$ which sends
a point $x=(z_0,\ldots,z_{m-1})$ to the point whose
coordinates are the values of the corresponding monomials
at $x$.

Let $Z\subset W_{p,q}\times \C^{m}$ be the set
\[( F_0, F_1,\ldots, F_n, x\ |\ F_i(x)=0 \text{\ for all\ } i).\]
There is a projection map $Z\to\Sigma$ which ``forgets'' the point $x$
where all the polynomials $F_i$ vanish.
Denote by $Z^{\Delta}$ the space of the simplicial resolution
associated to this map, with respect to the embedding
$Z\hookrightarrow W_{p,q}\times V$ which sends $(F_i, x)$
to $(F_i, v_{p,q}(x))$. We shall also use the corresponding non-degenerate
resolution, which will be denoted by $\widetilde{Z}^{\Delta}$.
The resolution $Z^{\Delta}$ will be, however, of greater importance to us
as it satisfies the inequality $(\ref{eq2})$ below.

Let $\widehat{Z}$, $\widehat{Z}^{\Delta}$ and $\widehat{\Sigma}$ be 
the one-point compactifications of $Z$, $Z^{\Delta}$ and $\Sigma$
respectively. The projection of $Z^{\Delta}$ onto $\Sigma$ extends to a
homotopy equivalence between $\widehat{Z}^{\Delta}$ and $\widehat{\Sigma}$. 
There is an increasing filtration  
\[\widehat{Z}_0\subset \widehat{Z}_1\subset\ldots\subset\widehat{Z}^{\Delta}\]
coming from the filtration $Z_r$ on $Z^{\Delta}$: the term $\widehat{Z}_0$ 
is the added point and $\widehat{Z}_r$ is equal to $Z_r\cup\widehat{Z}_0$.
In particular, for all $r>0$ the spaces $\widehat{Z}_r\backslash\widehat{Z}_{r-1}$
and $Z_r\backslash Z_{r-1}$ coincide.

The filtration on $\widehat{Z}^{\Delta}$ gives rise to a spectral sequence 
converging to the cohomology of
$\widehat{\Sigma}$:
\[E_1^{r,s}=H^{r+s}(\widehat{Z}_r,\widehat{Z}_{r-1},\Z),\]
where $\widehat{Z}_{-1}$ is the empty set.
The cohomology of $\widehat{\Sigma}$ is related  by the Alexander duality to 
the homology of $\Rt_f(p,q)$:
\[H^{r}(\widehat{\Sigma},{\bf Z})\simeq 
\widetilde{H}_{2N_{p,q}-r-1}(\Rt_f(p,q))\]
so the following spectral sequence converges to the reduced homology of $\Rt_f(p,q)$:
\begin{equation}\label{eq1}
E^1_{-r,s}=H^{2N_{p,q}+r-s-1}(\widehat{Z}_{r},\widehat{Z}_{r-1},\Z)
\end{equation}
with $r>0$, $s\geqslant 0$.

For $r>0$ the space $\widehat{Z}_r/\widehat{Z}_{r-1}=Z_r/Z_{r-1}$ is the one-point 
compactification of the space $Z_r\backslash Z_{r-1}$, which admits a rather explicit 
description, at least for small values of $r$.

A point of $Z_r$ is given by specifying (1) a point in $W_{p,q}$
which corresponds to $n+1$ polynomials $F_i$ vanishing simultaneously
on a non-empty subset of $\C^m$; (2) a point in $V$ which belongs to the
convex hull of at most $r$ points of the form $v_{p,q}(x_j)$
where $x_j$ are distinct points of $\C^m$ on which all the $F_i$ vanish.

The condition that a polynomial in $W^{i}_{p,q}$ vanishes at a
given point gives one linear inhomogeneous condition on its coefficients.
More generally, the condition that a $(p,q)$-polynomial in $W^{i}_{p,q}$
vanishes at $r$ distinct points $x_j$ produces exactly $r$ independent
conditions on its coefficients if and only if the convex hull of the points
$v_{p,q}(x_j)$  in $V$ is an $(r-1)$-dimensional simplex.

Hence,
\begin{equation}\label{eq2}
\begin{array}{cl}
\dim{Z_r\backslash Z_{r-1}}&
\leqslant 2(N_{p,q}-r(n+1))+2mr+(r-1)\\
&=2N_{p,q}-2r(n-m+1)+r-1.
\end{array}
\end{equation}

In general, little else can be said without a sophisticated analysis.
For example, a set of distinct points in $\C^m$ can give rise
to a set of linearly dependent conditions on the coefficients of a polynomial
vanishing at these points. However, the following is true:

\begin{prop}\label{Zr}
For all
$r\leqslant[\frac{p+1}{2}]$ the space $Z_r\backslash Z_{r-1}$ is
homeomorphic to a real vector bundle of rank
\[2N_{p,q}-(2n+1)r-1 \]
over the configuration space $C_{r}(\C^m)$  of $r$
distinct unordered points in $\C^m$.
\end{prop}

\noindent The proof of this statement is based on the following fact:

\begin{lem}\label{VdM}
If $r\leqslant p+1$ the images in $V$ under $v_{p,q}$ of
any set of $r$  distinct points in $\C^m$ span an
$(r-1)$-simplex.
\end{lem}
Indeed, by a linear change of coordinates in $\C^m$ one can acheive that the
values of the coordinate $z_0$ for all $r$ points are pairwise distinct;
then the Vandermonde matrix constructed of the powers of $z_0$ is
non-degenerate.

\begin{rem}\label{VdMrem}
Lemma~\ref{VdM} implies that for $r\leqslant p+1$ the condition that a
polynomial $F_i \in W^{i}_{p,q}$ vanishes at $r$ given points of $\C^m$
determines $r$ independent linear inhomogeneous conditions on the
coefficients of $F_i$. In general, these conditions may be incompatible.
The case $m=1$, $q=0$ is an example: a non-trivial polynomial of
degree $p$ in one variable cannot have $p+1$ roots. However, for
$r\leqslant p$ the $r$ affine hyperplanes in  $W^{i}_{p,q}$ defined by the
vanishing conditions at $r$ points are necessarily in general position.
This is proved by the same argument as Lemma~\ref{VdM}.
\end{rem}

%
\begin{proof}[Proof of Proposition~\ref{Zr}]
It follows from Lemma~\ref{VdM} that if $r\leqslant[\frac{p+1}{2}]$ then 
any two $(r-1)$-simplices in $V$ whose vertices are in the image of $v_{p,q}$ 
are either disjoint or have a common face. Hence, for these values of $r$,
given a point in the interior of the convex hull of $r$ distinct
points of the form $v_{p,q}(x_j)$ in $V$ one can determine the points $x_j$
up to order. Therefore, there exists a map
\[Z_r\backslash Z_{r-1}\to C_r(\C^m)\]
which keeps track of $r$-tuples of points in $\C^m$ on which the polynomials
$F_i$ vanish.

It was mentioned before that the condition for a polynomial to
vanish at a given point, determines an affine hyperplane in $W^{i}_{p,q}$.
By the remark to Lemma~\ref{VdM}, any set of $r$ points in $\C^m$
with $r\leqslant p$, determines a set of hyperplanes in general position, 
hence their intersection has complex codimension $r$ in $W^{i}_{p,q}$. 
Therefore, for all $r\leqslant [\frac{p+1}{2}]$ the fibre of the projection map
of $Z_r\backslash Z_{r-1}$
to  $C_r(\C^m)$ is a product of a complex vector space of complex dimension
$N_{p,q}-(n+1)r$, with an interior of an $(r-1)$-simplex. Verifying the local
triviality property is straightforward.
%
\end{proof}

The cohomology groups of the space $Z_r/Z_{r-1}$ are the same as those of
a reduced Thom space of a vector bundle of rank $2N_{p,q}-(2n+1)r-1$
over the one-point compactification $\widehat{C}_{r}(\C^m)$ of
$C_{r}(\C^m)$. In particular,
\[H^{i}(\widehat{C}_{r}(\C^m),{\bf Z})
=\widetilde{H}^{2N_{p,q}-(2n+1)r-1+i}(Z_r/Z_{r-1},{\bf Z})\]
for all $i$, and, hence,
\[E^1_{-r,s}=
H^{2(n+1)r-s}(\widehat{C}_{r}(\C^m),{\bf Z})\]
for all $0<r\leqslant[\frac{p+1}{2}]$ and all $s$. Note that this expression
does not depend on $p$ and $q$.

The term $E^1_{-r,s}$ of the spectral sequence for the homology
of $\Rt_f(p,q)$ is shown in the figure. It follows from (\ref{eq1})
and (\ref{eq2}) that all non-zero entries are
situated in the sector $-r < 0$, $s\geqslant 2(n-m+1)r$.
The entries in the strip $r\leqslant[\frac{p+1}{2}]$ are ``stable'' in the sense
that they do not depend on $p$ and $q$ and are preserved by the
stabilization maps. Let us make the last statement more precise.

Let us introduce the notation $\Sigma_{p,q}$ and $Z_{p,q}$ instead of
$\Sigma$ and $Z$, respectively. The inclusion map
\[\Rat_f(p,q)\to\Rat_f(p+1,q+1)\]
can be lifted to a map
\[\Rt_f(p,q)\to\Rt_f(p+1,q+1)\]
which multiplies all the $(n+1)$ $(p,q)$-polynomials by
$z_0\bar{z}_0+\ldots +z_m\bar{z}_m$. This map, in fact, extends to a
map $W_{p,q}\to W_{p+1,q+1}$; it sends
the discriminant $\Sigma_{p,q}$ to $\Sigma_{p+1,q+1}$. There is
also a corresponding map from $Z_{p,q}$ to $Z_{p+1,q+1}$. However, there is
no induced map between the degenerate simplical resolutions
$Z^{\Delta}_{p,q}$.

Consider the non-degenerate simplical resolution
$\widetilde{Z}^{\Delta}=\widetilde{Z}^{\Delta}_{p,q}$. There are maps
\begin{equation}\label{eq3}
Z^{\Delta}_{p,q}\stackrel{\pi_1}{\longleftarrow}\widetilde{Z}^{\Delta}_{p,q}
\stackrel{\zeta}{\to}
\widetilde{Z}^{\Delta}_{p+1,q+1}
\stackrel{\pi_2}{\longrightarrow}
{Z}^{\Delta}_{p+1,q+1}
\end{equation}
where $\pi_1$ and $\pi_2$ are homotopy equivalences
and the map $\zeta$ is induced by the stabilization map.

The first $[\frac{p+1}{2}]$ terms of the increasing filtration on
$\widetilde{Z}^{\Delta}_{p,q}$ coincide with the corresponding terms
for the degenerate resolution $Z^{\Delta}_{p,q}$. Hence, the $E^1$-terms
of the spectral
sequences associated with the four resolutions above all coincide
in the strip $r\leqslant[\frac{p+1}{2}]$. It is also obvious that the maps
$\pi_1$ and $\pi_2$ in (\ref{eq3}) induce the identity map in this strip.
The fact that  for $r\leqslant[\frac{p+1}{2}]$ the homomorphism of the $E^1$-terms
induced by $\zeta$ is identity as well, follows from the
explicit description of the filtration on $Z^{\Delta}_{p,q}$ obtained
above; essentially this is just the Thom isomorphism.

Due to the action of the differentials which
connect the stable and the unstable parts of the spectral sequence,
the stable region for the terms $E^2, E^3,\ldots$ is smaller than
that for the term $E^1$. A straightforward check shows that the set
\[
\begin{array}{l}
r\leqslant\left[\frac{p+1}{2}\right],\\
2(n-m+1)r\leqslant s \leqslant 
(2n-2m+1)\left(\left[\frac{p+1}{2}\right]+1\right)+r
\end{array}
\]
is in the stable region for the term $E^{\infty}$.

\begin{figure}
\[\epsffile{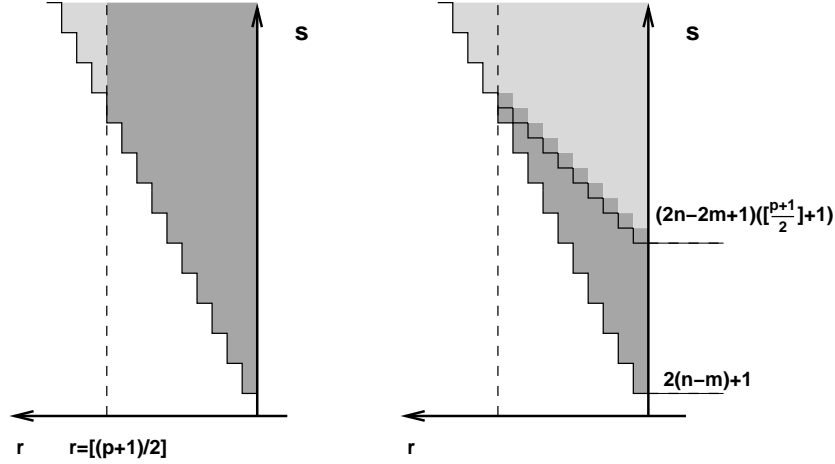}\]
\caption{The terms $E^{1}$ and $E^{\infty}$
of the spectral sequence converging to the  homology of 
$\Rat_f(p,q)$. The group $E^*_{a,b}$
is placed in the square \mbox{$a<x<a+1$}, \mbox{$b-1<y<b$}.
The shaded entries are stable.}
\end{figure}

For all $i<(2n-2m+1)([\frac{p+1}{2}]+1)$ the graded group
associated to $H_i(\Rat_f(p,q))$ is a direct sum of stable entries of the
term $E^{\infty}$ of the spectral sequence.
Hence, we obtain
\begin{prop}
The map
\[H_i(\Rat_f(p,q))\to H_i(\Rat_f(p+1,q+1))\]
induced by the stabilization map
is an isomorphism for all
\[ i<(2n-2m+1)\left(\left[\frac{p+1}{2}\right]+1\right).\]
\end{prop}

Together with Proposition~\ref{th:real} this implies Theorem~\ref{the}.


\section{Spaces of free maps.}\label{F}
The proof of Theorem~\ref{free} consists in applying Theorem~\ref{the} fibrewise
to the map
$$R:\Rat_{(d)}^{m,n}\to\Rat_{(d)}^{m-1,n}$$
induced by  restriction to a hyperplane. In what follows we denote by $D$
the stabilization dimension  
$(2n-2m+1)\left(\left[\frac{d+1}{2}\right]+1\right)$.

The restriction to a hyperplane 
$r:\text{\sl Map}_d(m,n)\to\text{\sl Map}_d(m-1,n)$ 
is a Serre fibration with the fibre $\Omega^{2m}\CPN{n}$. 
Denote by $t: T\to \Rat_{(d)}^{m-1,n}$ the pullback 
of $r$ to $\Rat_{(d)}^{m-1,n}$. 
A comparison of the spectral sequences for both 
fibrations shows that the inclusion $T\to\text{\sl Map}_d(m,n)$ induces 
isomorphisms in homology in dimensions smaller than $D$. Therefore, it is sufficient 
to prove that the inclusion $\Rat_{(d)}^{m,n}\to T$ is an isomorphism on 
homology in these dimensions.

According to Theorem~\ref{the}, the fibres of the map 
$R$ have the same homology as 
$\Omega^{2m}\CPN{n}$ in dimensions smaller than $D$. Moreover, 
we have the following:

\begin{lem}
There exists a finite filtration $F_0 \subset F_1 \subset\ldots
\subset F_q=\Rat_{(d)}^{m-1,n}$ 
such that  the map $R$ is a fibration over each connected 
component of $F_{i+1}\backslash F_i$. The spaces $F_{i}$ and $R^{-1}F_{i}$
can be triangulated so that each $F_{i-1}$ and each $R^{-1}F_{i-1}$ 
are subcomplexes in $F_{i}$ and $R^{-1}F_{i}$ respectively. 
\end{lem}

%
\begin{proof}
There exist smooth complex algebraic varieties $M$ and $N$ and a map
$R':M\to N$ such that  $\Rat_{(d)}^{m,n}$ and $\Rat_{(d)}^{m-1,n}$
are complements to closed subvarieties of $M$ and $N$ respectively, and such 
that $R'$ restricts to $R$ on  $\Rat_{(d)}^{m,n}$.

Indeed, the space $\Rat_{(d)}^{m,n}$ is a complement to a discriminant 
in a complex projective space;  the homogeneous coordinates in this projective 
space are 
the coefficients of the polynomials representing a point in $\Rat_{(d)}^{m,n}$.
From this point of view, the map $R$ is just a linear projection. Blowing up 
the 
indeterminacy locus of this projection one gets the variety $M$; the variety 
$N$ is 
just a complex projective space.  

There exist stratifications on $M$ and $N$, with $\Rat_{(d)}^{m,n}$ being the
maximal stratum of $M$, which turn $R':M\to N$ into a stratified map,
see Section~I.1.7 of \cite{GM}.   Take $F_{i}$ to be the filtration on 
$\Rat_{(d)}^{m-1,n}$ produced by intersecting $\Rat_{(d)}^{m-1,n}$ with the
above stratification of $N$. It follows from Thom's first isotopy lemma 
(see I.1.5 
of \cite{GM}) that $R$ is a fibration over each component of  
$F_{i+1}\backslash F_i$. The second part of the Lemma follows from the fact that 
every stratified set has a triangulation compatible with the stratification. 
%
\end{proof}

The inclusion  
$$R^{-1}(F_{i+1}\backslash F_{i}) \hookrightarrow 
t^{-1}(F_{i+1}\backslash F_{i})$$
is a morphism of fibrations over the same base. Hence, by 
Theorem~\ref{the},  it induces isomorphisms of homology
groups in dimensions less than  $D$. Consider the Mayer-Vietoris sequence
relating the homology of 
$R^{-1}F_{i+1}/R^{-1}F_{i}$ to 
the homology of $R^{-1}(F_{i+1}\backslash F_{i})$,
$R^{-1}V/ R^{-1}F_{i}$ and $R^{-1}\partial V$ where $V$ is
a small neighbourhood of $F_i$ in  $F_{i+1}$. (Notice 
that $R^{-1}V/ R^{-1}F_{i}$ is contractible and that 
$R:R^{-1}\partial V\to\partial V$ is a fibration.)
Comparing it
to the analogous Mayer-Vietoris sequence with $R$ replaced by $t$, one
verifies that the homomorphisms
$$H_{k}(R^{-1}F_{i+1}/R^{-1}F_{i})\to 
H_{k}(t^{-1}F_{i+1}/t^{-1}F_{i}) $$
are isomorphisms for $k<D$.
Now, to finish the proof of Theorem~\ref{free}, compare the spectral sequences 
associated to the filtrations $t^{-1}F_{i}$ on $T$, and $R^{-1}F_{i}$ 
on $\Rat_{(d)}^{m,n}$. 

\section{Final remarks.}

The estimate for the stabilization dimension in Theorem~\ref{the} in the case
$m=1$ is weaker than the estimate given by Segal's theorem. It is probable
that Lemma~\ref{VdM} on which our estimate is based, can be replaced by
a stronger statement. However, the methods of the present paper are
sufficient to recover Segal's theorem, at least for $n>1$. In the case
$m=1$, $q=0$ the simplicial resolution $Z^{\Delta}_{p,0}$ is, in fact,
non-degenerate, as a polynomial of degree $p$ in one variable
can have at most $p$ roots. Thus the explicit description of $Z_r/Z_{r-1}$
as a Thom space of a bundle over a configuration space is valid
not only for $r\leqslant[\frac{p+1}{2}]$, but for all $r\leqslant p$. This leads
to a better estimate for the stabilization dimension.

There is no doubt that Theorem~\ref{the} can be also strengthened by
replacing ``homology'' with ``homotopy'' in the case $m=n$. This would
require an argument similar to that used by Segal in \cite{Segal}.
Such an argument, however, is beyond the scope of the present paper.

\subsection*{Acknowledgments}
A great part of the motivation for this work came from a conversation with
V.\ Vassiliev who told me that simplicial resolutions could be applied to the 
study of the spaces $\Rat_{f}^{m,n}$. I am also grateful to M.\ Guest, 
S.\ P\'erez Esteva, A.\ Verjovsky and many others for useful discussions. 
Finally, I would like to thank Max-Planck-Institut f\"ur Mathematik, Bonn for 
hospitality during the stay at which a substantial part of the work presented 
here was done.







\begin{thebibliography}{99}
\bibitem{Boyer+} C.\ P.\ Boyer, J.\ C.\ Hurtubise \and R.\ J.\ Milgram,
Stability theorems for spaces of rational curves,
International Journal of Mathematics 12 (2001) 223--262.
%
\bibitem{GM} M.\ Goresky, R.\ MacPherson,
Stratified Morse Theory, Springer-Verlag, Berlin-Heidelberg-New York-Tokyo, 1988.
%
\bibitem{Gromov} M.\ Gromov, Oka's principle for holomorphic sections of
elliptic bundles, Journal of the American Mathematical Society  2 (1989) 
851--897.
%
\bibitem{Guest-A} M.\ Guest,  The topology of the space of rational
curves on a toric variety, Acta Mathematica 174 (1995) 119--145.
%
\bibitem{Guest-B} M.\ Guest, A.\ Kozlowski \and K.\ Yamaguchi,
Spaces of polynomials with roots of bounded multiplicity,
Fundamenta Mathematicae 161 (1999) 93--117.
%
\bibitem{Havlicek}
J.\ Havlicek,
On spaces of holomorphic maps from two copies of the Riemann sphere to
complex Grassmannians,
Real algebraic geometry and topology (East Lansing, MI, 1993), 83--116,
Contemporary  Mathematics 182,
American Mathematical Society, Providence, RI, 1995.
%
\bibitem{JM}
J.\ Mostovoy, Spaces of rational loops on a real projective space,
Transactions of the American Mathematical Society 353 (2001) 1959--1970.
%
\bibitem{Kallel+}
S.\ Kallel \and R.\ J.\ Milgram,  The geometry of the space of holomorphic 
maps from a
Riemann surface to a complex projective space,
Journal of Differential Geometry  47 (1997) 321--375.
%
\bibitem{KYlin}
A.\ Kozlowski\ and\ K.\ Yamaguchi,  Spaces of holomorphic maps between complex
projective spaces of degree one,  Topology and its Applications 132 (2003)
139--145.
%
\bibitem{Rudin} W.\ Rudin,  Principles of mathematical analysis,
third edition, McGraw-Hill Book Co., New York-Auckland-D\"usseldorf, 1976.
%
\bibitem{Segal} G.\ Segal, The topology of spaces of rational
functions, Acta Mathematica  143 (1979), 39-72.
%
\bibitem{V1} V.\ A.\ Vassiliev,
Complements of discriminants of smooth maps: topology and applications,
revised edition, Translations of Mathematical Monographs 98,
American Mathematical Society, Providence, RI, 1994.
%
\bibitem{V2} {\cy V.\ A.\ Vasilp1ev},
{\cy Topologiya dopolnenii0 k diskriminantam,  Izdatelp1stvo
FAZIS, Moskva, 1997}.

\end{thebibliography}
\end{document}